\documentclass[smallextended,referee,envcountsect,]{svjour3}
\smartqed

\usepackage{graphicx}
\usepackage{amsmath,amssymb,amssymb,dsfont}
\usepackage{multirow,booktabs,tabularx,cases,threeparttable}
\usepackage{algorithm, algorithmic}
\usepackage[latin1]{inputenc}
\usepackage{bbding}
\usepackage[misc]{ifsym}
\usepackage[colorlinks,
           citecolor=blue,
           urlcolor=blue,
            linkcolor=blue,
            anchorcolor=blue]{hyperref}

\newcommand{\cl}{\textcolor{black}}

\journalname{JOTA}

\begin{document}

\title{A tensor optimization algorithm for computing Lagrangians of hypergraphs}

\author{Jingya Chang\textsuperscript{1}  \and
        Bin Xiao\textsuperscript{1}   \and
        Xin Zhang\textsuperscript{2}
}
\institute{
 \Letter   \quad  Xin Zhang   \\
          \phantom{four} zhangxin0619@126.com   \\
                  \\
          \phantom{four}  Jingya Chang  \\
           \phantom{four} jychang@gdut.edu.cn \\
           \\
           \phantom{four} Bin Xiao   \\
           \phantom{four} 2112114085@mail2.gdut.edu.cn    \\
           \\
          \textsuperscript{1}  \phantom{on}School of  Mathematics and Statistics, Guangdong University of Technology, \\
                \phantom{four}   Guangzhou 510006, People's Republic of China \\
           \\
            \textsuperscript{2} \phantom{on}School of Arts and Science, Suqian University, Suqian 223800, \\
            \phantom{four} People's Republic of China  \\
}

\date{Received: date / Accepted: date}

\maketitle

\begin{abstract}
The Lagrangian of a hypergraph is a crucial tool for studying hypergraph extremal problems. Though Lagrangians of some special structure hypergraphs have closed-form solutions, it is a challenging problem to compute the Lagrangian of a general large scale hypergraph. In this paper, we exploit a fast computational scheme involving the adjacency tensor of a hypergraph. Furthermore, we propose to utilize the gradient projection method \cl{on a simplex} from nonlinear optimization for solving the Lagrangian of a large scale hypergraph iteratively. Using the {\L}ojasiewicz gradient inequality, we analyze the global and local convergence of the gradient projection method. Numerical experiments illustrate that the proposed numerical method could compute Lagrangians of large scale hypergraphs efficiently.
\end{abstract}

\newcommand{\av}{\mathbf{a}}
\newcommand{\bv}{\mathbf{b}}
\newcommand{\gv}{\mathbf{g}}
\newcommand{\vv}{\mathbf{v}}
\newcommand{\xv}{\mathbf{x}}
\newcommand{\yv}{\mathbf{y}}
\newcommand{\At}{\mathcal{A}}
\newcommand{\Pt}{\mathcal{P}}
\newcommand{\one}{\mathbf{e}}
\newenvironment{funding}{\begin{funding}}
{\end{funding}}

\newtheorem{Thm}{Theorem}[section]
\newtheorem{Def}[Thm]{Definition}
\newtheorem{Ass}[Thm]{Assumption}
\newtheorem{Lem}[Thm]{Lemma}
\newtheorem{Cor}[Thm]{Corollary}

\section{Introduction}

Hypergraphs have important applications in science and engineering, such as subspace clustering \cite{BP13,CQZ17,ZLZ22}, hypergraph matching \cite{HPY22}, image processing \cite{CCQY20}, and network analysis \cite{BGL16,LPL22}, due to its capability of modeling multiwise similarity. Hypergraph extremal problems, which maximize values of constrained multilinear functions, attract attentions of scholars from many fields of mathematics \cite{GLM-20}. In this paper, we focus on one kind of hypergraph extremal problems: the Lagrangian of a hypergraph.

\begin{figure}[h]
  \centering
  \includegraphics[width=5cm]{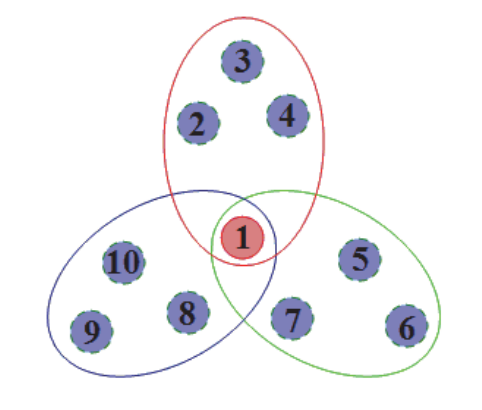}\\
  \caption{A 4-uniform hypergraph: sunflower}\label{sunflower}
\end{figure}

An $r$-uniform hypergraph ($r$-graph) $G=(V,E)$ consists of a vertex set $V:=\{1,2,\dots,n\}$ and an edge set $E\subseteq V^{(r)}$, where $V^{(r)}$ denotes the collection of all subsets of $V$ of size $r$,
$m=|E|$ is the number of edges. For example the sunflower hypergraph in Figure \ref{sunflower} is a $4$-uniform hypergraph with $10$ vertices. The weight polynomial for $G$ is defined as $w(G,\xv):=\sum_{e\in E}\prod_{i\in e}x_i$, where legal weighting for $G$ satisfies $x_i\ge0$ for all $i\in V$ and $\sum_{i\in V}x_i=1$.
Let $\Delta:=\{\xv\in \Re^n_+:\one^T\xv=1\}$ be the legal weighting set that is a simplex. Here $\one$ is the all one vector. The Lagrangian of an $r$-graph $G$ \cite{MS-65,Ta-02} is to maximize the weight polynomial $w(G,\xv)$ under the simplex constraint
\begin{equation}\label{Lagr}
  \lambda(G) := \max\left\{w(G,\xv)=\sum_{e\in E}\prod_{i\in e}x_i~:~\one^T\xv=1,~\xv\in\Re^n_+\right\}.
\end{equation}

Motzkin and Straus \cite{MS-65} introduced the concept of Lagrangians for 2-graphs and built a bridge between Lagrangians and Tur\'{a}n's theorem. For $r=2$, $\lambda(G)$ is achieved by equally distributing the weight over the vertices of the largest clique in $G$ and setting $x_i=0$ for all other vertices $i$. Moreover, Lagrangians of hypergraphs have some significant results.
Let $K_n^r$ denote the complete $r$-graph on $V:=\{1,2,\dots,n\}$, where $E=V^{(r)}$. Then, it is pointed out in \cite{HPW20} that
\begin{equation*}
  \lambda(K_n^r)={n \choose r}\frac{1}{n^r}
\end{equation*}
and the associated solution is $\xv^*=(\frac{1}{n},\dots,\frac{1}{n})^T$.
Let $t$ be a positive integer and $\{f_1, f_2, \dots, f_{t+2}\}$ be pairwise disjoint $3$-subsets. Let $\mathcal{F}^*$ be the $3$-graph with the vertex set $V = \bigcup_{1\le k\le t+2}f_k$ and the edge set $\{f_1, f_2, \dots, f_{t+2}\}\cup\{f\in V^3 : |f \cap f_k|\le 1\text{ for }k=1,2,\dots,t+2\}$. Then,
\begin{equation*}
  \lambda(\mathcal{F}^*)=\left(27{t+2 \choose 3}+t+2\right)\frac{1}{27(t+2)^3}<\lambda(K_{3t+4}^3).
\end{equation*}
Frankl and F\"{u}redi \cite{FF-89} conjectured in 1989 that an initial segment of the colexicographic order has the largest Lagrangian of any $r$-graph with size $m$. The Frankl--F\"{u}redi conjecture was studied by many authors \cite{LZ15,STZP14}. Although the Frankl--F\"{u}redi conjecture is true in the case of $r=3$ \cite{Ta-02,GLM-21}, Gruslys et al.  \cite{GLM-20} disproved it for $r\ge4$.

While theoretical researches on Lagrangians of hypergraphs are rich, the values of Lagrangians of ordinary hypergraphs are still unclear with the aid of theoretical analysis. Therefore, we are going to design numerical method for computing the Lagrangian of general uniform hypergraphs.

Rich theoretical results imply that the Lagrangian of a hypergraph reveals special nature of the adjacency tensor of a hypergraph. However, it is a challenging problem to compute the Lagrangian of a large scale general hypergraph. In this paper, we customize the gradient projection algorithm from nonlinear optimization for computing the Lagrangian of a large scale general hypergraph. On one hand, from the viewpoint of spectral hypergraph theory, the Lagrangian of an $r$-graph $G$ is closely related to the adjacency tensor $\At$ of $G$, which is an $r$th order structure tensor with valuable symmetry and sparsity. Fast computations of the structure tensors arising from a hypergraph could be employed for the Lagrangian problem. The computation cost for computing function values and gradients of the weight polynomial $w(G,\xv)$ is proportional to the size $m$ of the $r$-graph and the square of $r$. On the other hand, the legal weighting set $\Delta$ is a simplex, which is a closed convex set. Since the projection onto the simplex is cheap, we design a gradient projection method for solving the Lagrangian of an $r$-graph $G$, where the initial step size at each iteration is Barzilai--Borwein step \cl{size}.

Due to the semi-algebraic property of the Lagrangian of an $r$-graph $G$, the {\L}ojasiewicz inequality holds. Using the {\L}ojasiewicz inequality, we analyze the global convergence of the gradient projection algorithm, i.e., the sequence of iterates generated by the gradient projection algorithm converges to a critical point with linear or sublinear rate. Numerical experiments on small and large scale hypergraphs illustrate that the gradient projection algorithm is powerful and efficient. In particular, the gradient projection algorithm could compute Lagrangians of hypergraphs with thousands of edges.

The outline of this paper is drawn as follows. Section 2 presents the gradient projection algorithm and associated fast computations on hypergraphs. Global and local convergence is analyzed in Section 3. Numerical experiments on small and large scale hypergraphs are reported in Section 4. Finally, some concluding remarks are made in Section 5.

\section{Gradient projection method}

To handle large scale hypergraphs, there are roughly two kinds of strategies: (i) hardware acceleration develops and optimizes CPU/GPU kernels to process hypergraph algorithms and (ii) software acceleration uses tensors and tensor operators to represent hypergraph computing into a unique (compact) format that can be executed efficiently \cite{K21}. In this paper, we follow the software acceleration strategy to exploit tensor representations and fast computations for the purpose of the Lagrangian computing of hypergraphs. At the beginning, for an $r$-graph $G$, the weight polynomial $w(G,\xv)$ is determined by the adjacency tensor of $G$. First we introduce the definition of tensor and hypergraph related adjacency tensor.
\begin{Def}[Symmetry tensor \cite{QL-book}]
  A tensor $$\mathcal{T}=(t_{i_1\ldots i_r})\in \Re^{[r,n]}, \qquad \text{for} \ i_j = 1,\ldots, n, j = 1,\ldots,r$$ is an $r$th order $n$ dimensional symmetric tensor if the value of $t_{i_1\ldots i_k}$ is invariable under any permutation of its indices.
\end{Def}

\begin{Def}[Adjacency tensor \cite{CoD-12}]
  Let $G=(V,E)$ be an $r$-graph with $n$ vertices. The adjacency tensor of $G$ is an $r$th order $n$-dimensional symmetric tensor $\At=[a_{i_1 \cdots i_r}]$, of which elements are
  \begin{equation*}
    a_{i_1 \cdots i_r}=\left\{\begin{aligned}
      &\frac{1}{(r-1)!} && \quad\text{ if }\{i_1,\ldots,i_r\}\in E, \\
      &0                && \quad\text{ otherwise. }
    \end{aligned}\right.
  \end{equation*}
\end{Def}

Utilizing the adjacency tensor $\At$ of the $r$-graph $G$, we have the following lemma.

\begin{Lem}
  Let $G$ be an $r$-graph with  $n$ vertices. Its weight polynomial $w(G,\xv)$ could be rewritten as
  \begin{equation*}
    w(G,\xv) = \frac{1}{r}\At\xv^r,
  \end{equation*}
  where $\At\xv^r := \sum_{i_1=1}^n\cdots\sum_{i_r=1}^n a_{i_1\cdots i_r}x_{i_1}\cdots x_{i_r}$.
\end{Lem}
\begin{proof}
  By direct calculations, it yields that
  \begin{eqnarray}
  \frac{1}{r}\At\xv^r &=& \frac{1}{r}\sum_{i_1=1}^n\cdots\sum_{i_r=1}^n a_{i_1\cdots i_r}x_{i_1}\cdots x_{i_r} \nonumber\\
    &=& \frac{1}{r}\sum_{{\{i_1,\dots,i_r\}\in E}\atop i_1<\dots<i_r} \frac{1}{(r-1)!} r! x_{i_1}\cdots x_{i_r} \nonumber\\
    &=& \sum_{e\in E}\prod_{i\in e}x_i = w(G,\xv). \label{aaabbb}
  \end{eqnarray}
  The proof is completed.
\end{proof}

Hence, to compute the Lagrangian of an $r$-graph $G$, we will solve the following optimization problem
\begin{equation}\label{OPM}
\lambda(G)=\left\{\begin{aligned}
  & \max && f(\xv)=\frac{1}{r}\At\xv^r \\
  & ~\mathrm{s.t.}&& \xv\in\Delta,
\end{aligned}\right.
\end{equation}
which is equivalent to \eqref{Lagr}. Since the adjacency tensor $\At$ of $G$ is symmetric, it holds \cite{QL-book,CDQY18} that $\nabla f(\xv)=\At\xv^{r-1}$, where elements of the vector $\At\xv^{r-1}$ are
\begin{equation*}
  (\At\xv^{r-1})_j=\sum_{i_2=1}^n\cdots\sum_{i_r=1}^n a_{ji_2\cdots i_r}x_{i_2}\cdots x_{i_r} \qquad
  \forall~j=1,\dots,n.
\end{equation*}

\begin{figure}[t]
  \centering
  \includegraphics[width=0.6\textwidth]{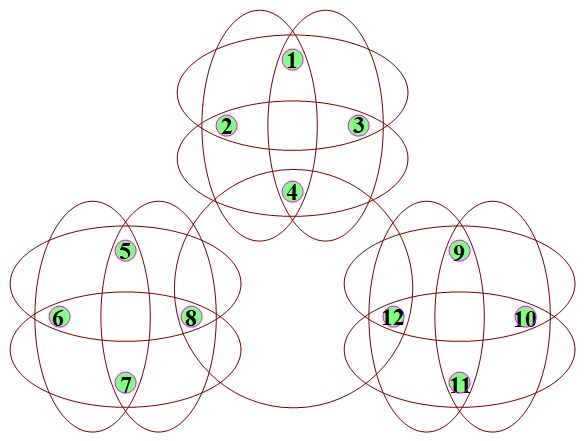}\\
  \caption{A toy hypergraph.}\label{NE:toy}
\end{figure}

Since the adjacency tensor $\At$ of an $r$-graph $G$ is a structure tensor, the storage and computation involving $\At$ are economic \cite{CCQ16}. For example, we consider the $3$-graph illustrated in Figure \ref{NE:toy}, i.e., $r=3$ here. The $3$-graph has $n=12$ vertices and $m=13$ edges. The associated adjacency tensor is a $3$rd order $12$ dimensional symmetric tensor with $n^r=1,728$ elements. Using the storage technique from \cite{CCQ16}, we save the incidence matrix of the $3$-graph in a compact form
\begin{equation*}
  H=\left(\begin{array}{ccccccccccccc}
        1 & 1 & 1 & 2 & 5 & 5 & 5 & 6 &  9 &  9 &  9 & 10 &  4 \\
        2 & 2 & 3 & 3 & 6 & 6 & 7 & 7 & 10 & 10 & 11 & 11 &  8 \\
        3 & 4 & 4 & 4 & 7 & 8 & 8 & 8 & 11 & 12 & 12 & 12 & 12 \\
    \end{array}\right)^T,
\end{equation*}
which only has $m\times r=39$ positive integers.

To compute the scalar $\At\xv^r$ for an input vector $\xv=[x_j]\in\Re^n$, we introduce a matrix $Y=[y_{ik}]\in\Re^{m\times r}$ of which elements are
\begin{equation*}
  y_{ik} = x_{H_{ik}}, \qquad\forall~i=1,\dots,m\text{ and }k=1,\dots,r.
\end{equation*}
Then, it can be deduced from \eqref{aaabbb} that
\begin{equation*}
  \At\xv^r = r\sum_{i=1}^m\prod_{k=1}^r y_{ik}.
\end{equation*}
Similarly, to compute the vector $\At\xv^{r-1}$, we define
\begin{equation*}
  M^{(k)} := [\delta(i,H_{jk})]_{ij} \quad\text{ and }\quad
  \yv^{(k)} := \left[\prod_{\ell=1\atop \ell\neq k}^r y_{j\ell}\right]_j,
\end{equation*}
for $k=1,\dots,r$, where $\yv^{(k)}$ specifies values and $M^{(k)}$ indicates locations. Here, $\delta(\cdot,\cdot)$ stands for the Kronecker delta. Then, we have
\begin{equation*}
  \At\xv^{r-1} = \sum_{k=1}^r M^{(k)}\yv^{(k)}.
\end{equation*}
The computational cost for computing $\At\xv^r$ and $\At\xv^{r-1}$ is about $\mathcal{O}(mr^2)$, which is cheap due to the structure of an $r$-graph.

It is easy to see that the feasible region $\Delta$ is a closed convex set. The Euclidean projection
\begin{equation*}
  \pi_{\Delta}(\av) := \arg\min\limits_{\xv\in\Delta} \|\xv-\av\|
\end{equation*}
exists and is unique for any vector $\av\in\Re^n$. Algorithm \ref{Alg:Proj} provides an $\mathcal{O}(n\log n)$ approach for computing the Euclidean projection onto $\Delta$ \cite{WC'13,CY'11}.

\begin{algorithm}[t]
\caption{Euclidean projection $\pi_{\Delta}(\av)$.}\label{Alg:Proj}
\begin{algorithmic}[1]
  \STATE Sort the input vector $\av\in\Re^n$ into an auxiliary vector $\bv$ such that $b_1\ge b_2\ge \dots\ge b_n$.
  \STATE Find the largest index $\ell=\max\{j: jb_j+1-\sum_{i=1}^j b_i>0, 1\le j\le n\}$.
  \STATE Set $\lambda = (1-\sum_{i=1}^{\ell} b_i)/{\ell}$.
  \STATE Output $\xv$ with $x_i=\max\{a_i+\lambda,0\}$.
\end{algorithmic}
\end{algorithm}

Now, we give the gradient projection algorithm formally in Algorithm \ref{Alg:pGrad}. The gradient projection algorithm is an iterative method from nonlinear optimization \cite{Ber-book}. Starting from an initial  iteration point $\xv_0\in\Delta$, we compute the gradient of the objective function and the associated gradient projection direction. Next, we select an initial step size by Barzilai--Borwein method \cite{BB88,HDL21} and then perform backtracking line search along this gradient projection direction. At each iteration, only one projection is performed. If $\gv_c$ in Step 3 of Algorithm \ref{Alg:pGrad} is sufficiently small or the total number of iterations is large enough, we terminate the algorithm.

\begin{algorithm}[t]
\caption{A gradient projection algorithm (GPA).}\label{Alg:pGrad}
\begin{algorithmic}[1]
  \STATE Choose $0<\underline{\alpha}\le\alpha_0\le\overline{\alpha}$, $\beta\in(0,1)$, $\eta\in(0,1/2]$, and $\xv_0\in\Delta$. Set $c\gets0$.
  \WHILE{the sequence of iterates does not converge}
    \STATE Set $\gv_c = \pi_{\Delta}(\xv_c + \alpha_c\nabla f(\xv_c))-\xv_c$.
    \STATE Find the smallest nonnegative integer $j$ such that
      \begin{equation}\label{LineSch}
        f(\xv_c+\beta^j\gv_c)-f(\xv_c) \ge \eta\beta^j \gv_c^T \nabla f(\xv_c).
      \end{equation}
    \STATE Set $\rho_c=\beta^j$ and $\xv_{c+1}=\xv_c+\rho_c\gv_c$.
    \STATE Choose $\alpha_{c+1}\in[\underline{\alpha},\overline{\alpha}]$ and set $c\gets c+1$.
  \ENDWHILE
\end{algorithmic}
\end{algorithm}

\section{Convergence analysis}


Now, we analyze the convergence of the proposed GPA. Part of our analysis on basic theory is similar to \cite{Ber-book,HZ'06} and we present here for the purpose of completion. First, we call $\xv_*\in\Delta$ a critical point of \eqref{OPM} if
\begin{equation*}
  \pi_{\Delta}(\xv_*+\nabla f(\xv_*))=\xv_*.
\end{equation*}

Since the simplex $\Delta$ is a nonempty, closed, and convex set, we get the following lemma \cite[Proposition 2.1]{HZ'06}.

\begin{Lem}\label{Lem:4-1}
  Let $\gv^{\alpha}(\xv)=\pi_{\Delta}(\xv+\alpha\nabla f(\xv))-\xv$. Then, \\
  1) $\|\gv^{\alpha}(\xv)\|$ is nondecreasing in $\alpha>0$ for all $\xv\in\Delta$; \\
  2) $\|\gv^{\alpha}(\xv)\|/\alpha$ is nonincreasing in $\alpha>0$ for all $\xv\in\Delta$; \\
  3) $\gv^{\alpha}(\xv)^T\nabla f(\xv)\ge\|\gv^{\alpha}(\xv)\|^2/\alpha$ for all $\xv\in\Delta$ and $\alpha>0$; \\
  4) for any $\xv\in\Delta$ and $\alpha>0$, $\gv^{\alpha}(\xv)=0$ if and only if $\xv$ is a \cl{stationary} point for \eqref{OPM}.
\end{Lem}

Because $f(\xv)=\frac{1}{r}\At\xv^r$ is a polynomial and the simplex $\Delta$ is compact, the Hessian $\nabla^2 f(\xv)$ is bounded, i.e., there exists a constant $M\ge1$ such that
\begin{equation*}
  \|\nabla^2 f(\xv)\| \le M  \qquad \forall \xv\in\Delta.
\end{equation*}

\begin{Lem}\label{Lem:4-2}
  There exists a constant $\kappa_1>0$ such that for any $\rho_c$ generated by Algorithm \ref{Alg:pGrad}
  \begin{equation*}
    \rho_c\ge \kappa_1.
  \end{equation*}

\end{Lem}
\begin{proof}
  We will prove that the inequality \eqref{LineSch} is valid
  when $\rho\in[0, \frac{2(1-\eta)}{\overline{\alpha}M}]$. By Lemma \ref{Lem:4-1}, it yields that $\|\gv_c\|^2\le\alpha_c\gv_c^T\nabla f(\xv_c)$. From Taylor's formula, we have
  \begin{eqnarray*}
    f(\xv_c+\rho\gv_c) &\ge& f(\xv_c)+\rho \gv_c^T\nabla f(\xv_c)-\frac{M}{2}\rho^2\|\gv_c\|^2 \\
      &\ge& f(\xv_c)+\rho \gv_c^T\nabla f(\xv_c) -\frac{M}{2}\rho^2\alpha_c\gv_c^T\nabla f(\xv_c) \\
      &\ge& f(\xv_c)+\rho \gv_c^T\nabla f(\xv_c) -(1-\eta)\rho\gv_c^T\nabla f(\xv_c) \\
      &=& f(\xv_c)+\eta\rho \gv_c^T\nabla f(\xv_c),
  \end{eqnarray*}
  where the last inequality holds owing to the assumption $0\le\rho\alpha_cM\le \rho\overline{\alpha}M\le 2(1-\eta)$. According to Step 4 of Algorithm GPA,  this lemma is valid if $\kappa_1:=\frac{2\beta(1-\eta)}{\overline{\alpha}M}$.
\end{proof}

\begin{Thm}\label{Thm:4-3}
  Suppose that $\{\xv_c\}$ is an infinity sequence of iterates generated by GPA. Then, we have
  \begin{equation*}
    \lim_{c\to\infty} \|\pi_{\Delta}(\xv_c + \nabla f(\xv_c))-\xv_c\| = 0.
  \end{equation*}
  That is to say, every limit point of $\{\xv_c\}$ is a critical point.
\end{Thm}
\begin{proof}
  From \eqref{LineSch}, Lemmas \ref{Lem:4-1} and \ref{Lem:4-2}, we obtain
  \begin{equation}\label{fun-lowbd}
    f(\xv_{c+1})-f(\xv_c) \ge \eta\rho_c\gv_c^T\nabla f(\xv_c) \ge \frac{\eta\rho_c}{\alpha_c}\|\gv_c\|^2
      \ge \frac{\eta\kappa_1}{\overline{\alpha}}\|\gv_c\|^2.
  \end{equation}
  Hence, we have
  \begin{equation*}
    \sum_{c=1}^{\infty} \|\gv_c\|^2 \le \frac{\overline{\alpha}}{\eta\kappa_1}\sum_{c=1}^{\infty}[f(\xv_{c+1})-f(\xv_c)]
       \le \frac{\overline{\alpha}}{\eta\kappa_1}\lambda(G),
  \end{equation*}
  which means that $\|\gv_c\| \to0$ as $c\to\infty$.

  On one hand, if $\alpha_c\ge1$, we get $\|\gv^{\alpha_c}(\xv_c)\|\ge\|\gv^1(\xv_c)\|$ by Lemma \ref{Lem:4-1}. On the other hand, $\alpha_c<1$ and hence $\|\gv^{\alpha_c}(\xv_c)\|/\alpha_c\ge\|\gv^1(\xv_c)\|$. Hence, we have
  \begin{equation}\label{gc-lowbd}
    \|\gv_c\|=\|\gv^{\alpha_c}(\xv_c)\|\ge \min\{\alpha_c,1\}\|\gv^1(\xv_c)\| \ge \min\{\underline{\alpha},1\}\|\gv^1(\xv_c)\|.
  \end{equation}
 Therefore $\|\gv^1(\xv_c)\|$ tends to zero. The proof is then completed.
\end{proof}

\subsection{Further results based on {\L}ojasiewicz gradient inequality}

Because the objective function and the constraint set of \eqref{OPM} are semi-algebraic, the following {\L}ojasiewicz gradient inequality is valid \cite{Loj-63,BDL-07}. The analysis of this subsection is based on work in \cite{ABRS10}.

\begin{Thm}[{\L}ojasiewicz property]\label{Thm:5-1}
  Suppose that $\xv_*$ is a critical point of \eqref{OPM}, i.e., $\gv^1(\xv_*)=0$. Then, there exist a neighborhood of $\xv_*$ denoted as $\mathds{U}(\xv_*)$, an exponent $\theta\in[1/2,1)$, and a positive constant $C$ such that the following inequality
  \begin{equation*}
    \|\gv^1(\xv)\| \ge C|f(\xv)-f(\xv_*)|^{\theta}
  \end{equation*}
  holds for all $\xv\in \mathds{U}(\xv_*) \cap\Delta$.
\end{Thm}

\begin{Lem}\label{Lem:5-2}
  Let $\xv_*$ be a limiting point of $\{\xv_c\}$. The initial iterate $\xv_0$ is sufficiently close to $\xv_*$ in the sense that $\xv_0\in\mathds{B}(\xv_*,\sigma)\cap\Delta$, where $\mathds{B}(\xv_*,\sigma):=\{\xv\in\Re^n:\|\xv-\xv_*\|<\sigma\}\subseteq\mathds{U}(\xv_*)$ and $\sigma\ge \frac{\max\{1,\overline{\alpha}\}}{\eta C(1-\theta)}|f(\xv_0)-f(\xv_*)|^{1-\theta}+\|\xv_0-\xv_*\|$. Then, we have the following two assertions:
  \begin{equation*}
    \xv_c\in\mathds{B}(\xv_*,\sigma)\cap\Delta \qquad\text{ for }c=0,1,2,\dots,
  \end{equation*}
  and
  \begin{equation*}
    \sum_{c=0}^{\infty}\|\xv_c-\xv_{c+1}\| \le \frac{\max\{1,\overline{\alpha}\}}{\eta C(1-\theta)}|f(\xv_0)-f(\xv_*)|^{1-\theta}.
  \end{equation*}
\end{Lem}
\begin{proof}
The conclusions can be  proved by  induction. It is obvious to see that $\xv_0\in\mathds{B}(\xv_*,\sigma)\cap\Delta$. Next, by supposing that there exists a positive integer $\ell$ such that
  \begin{equation*}
    \xv_c\in\mathds{B}(\xv_*,\sigma)\cap\Delta \qquad\text{ for }c=0,1,\dots,\ell,
  \end{equation*}
  we show $\xv_{\ell+1}\in\mathds{B}(\xv_*,\sigma)\cap\Delta$.

  Define $\phi(t):=\frac{1}{C(1-\theta)}|t-f(\xv_*)|^{1-\theta}$. It is easy to see that $\phi(t)$ is a monotonically decreasing and  \cl{concave} function for  $t<f(\xv_*)$. Then, for $0\le c\le \ell$, we have
  \begin{eqnarray*}
    \phi(f(\xv_c))-\phi(f(\xv_{c+1}))
      &\ge& \phi'(f(\xv_c))(f(\xv_c)-f(\xv_{c+1})) \\
      &=& \frac{-1}{C|f(\xv_c)-f(\xv_*)|^{\theta}} (f(\xv_c)-f(\xv_{c+1})) \\
      &\ge& \frac{1}{\|\gv^1(\xv_c)\|}(f(\xv_{c+1})-f(\xv_c)) \\
      &\ge& \frac{1}{\|\gv^1(\xv_c)\|}\frac{\eta\rho_c}{\alpha_c}\|\gv_c\|^2 \\
      &\ge& \frac{1}{\|\gv^1(\xv_c)\|}\frac{\eta\rho_c}{\alpha_c}\|\gv_c\|\min\{\alpha_c,1\}\|\gv^1(\xv_c)\| \\
      &=& \frac{\eta\min\{\alpha_c,1\}}{\alpha_c}\|\rho_c\gv_c\| \\
      &=& \eta\min\{1,\alpha_c^{-1}\}\|\xv_c-\xv_{c+1}\|,
  \end{eqnarray*}
  where the second inequality is obtained based on Theorem \ref{Thm:5-1}, the third inequality is deduced from \eqref{fun-lowbd}, and the fourth inequality is valid because of \eqref{gc-lowbd}. The above inequalities indicate that
  \begin{equation*}
    \|\xv_c-\xv_{c+1}\| \le \frac{\max\{1,\overline{\alpha}\}}{\eta}\left(\phi(f(\xv_c))-\phi(f(\xv_{c+1}))\right),
  \end{equation*}
  for $0\le c\le \ell$. Hence, it holds that
  \begin{equation}\label{iterate-length}
    \sum_{c=0}^{\ell} \|\xv_c-\xv_{c+1}\|
      \le \frac{\max\{1,\overline{\alpha}\}}{\eta}\sum_{c=0}^{\ell}\left(\phi(f(\xv_c))-\phi(f(\xv_{c+1}))\right)
      \le \frac{\max\{1,\overline{\alpha}\}}{\eta}\phi(f(\xv_0)).
  \end{equation}
  Therefore, we have
  \begin{equation*}
    \|\xv_{\ell+1}-\xv_*\| \le \|\xv_0-\xv_*\|+\sum_{c=0}^{\ell} \|\xv_c-\xv_{c+1}\|
      \le \|\xv_0-\xv_*\|+\frac{\max\{1,\overline{\alpha}\}}{\eta}\phi(f(\xv_0)) <\sigma,
  \end{equation*}
  which implies that $\xv_{\ell+1}\in\mathds{B}(\xv_*,\sigma)$. On the other hand, every iterate is feasible by the mechanism of GPA.  Hence, we obtain the first assertion. The second  conclusion comes from \eqref{iterate-length} straightforwardly by setting $\ell\to\infty$.
\end{proof}

\begin{Thm}\label{Thm:5-3}
  Suppose that GPA generates an infinite sequence of iterates $\{\xv_c\}$. Then the whole sequence $\{\xv_c\}$ converges to a critical point $\xv_*$.
\end{Thm}
\begin{proof}
  Owing to the compactness of the feasible region $\Delta$, there exists at least a limiting point $\xv_*$ of $\{\xv_c\}$. From Theorem \ref{Thm:4-3}, $\xv_*$ is a critical point of the optimization problem \eqref{OPM}. On the other hand, there is an iteration $c_0$ such that $\xv_{c_0}\in\mathds{B}(\xv_*,\sigma)\cap\Delta$, where $\sigma$ is specified by Lemma \ref{Lem:5-2}. By regarding $\xv_{c_0}$ as an initial iterate in GPA, we know
  \begin{eqnarray*}
    \sum_{c=0}^{\infty}\|\xv_c-\xv_{c+1}\| &=& \sum_{c=0}^{c_0-1}\|\xv_c-\xv_{c+1}\| + \sum_{c=c_0}^{\infty}\|\xv_c-\xv_{c+1}\| \\
      &\le& \sum_{c=0}^{c_0-1}\|\xv_c-\xv_{c+1}\| + \frac{\max\{1,\overline{\alpha}\}}{\eta C(1-\theta)}|f(\xv_{c_0})-f(\xv_*)|^{1-\theta}
      <\infty,
  \end{eqnarray*}
  where the first inequality is owing to the second assertion of Lemma \ref{Lem:5-2}. Hence, we claim that the whole sequence $\{\xv_c\}$ converges.
\end{proof}

\subsection{Convergence rate}

In this subsection, we analyze the convergence rate of our GPA algorithm based on works in \cite{AtB-09,HL'18}.

\begin{Thm}\label{Thm:5-4}
  Suppose GPA generates an infinite sequence of iteration points $\{\xv_c\}$ that converges to a critical point $\xv_*$. Then, we have the following estimations on convergence rate.
  \begin{itemize}
    \item If $\theta=1/2$, there exist $\gamma>0$ and $\mu\in(0,1)$ such that
      \begin{equation*}
        \|\xv_c-\xv_*\| \le \gamma\mu^c.
      \end{equation*}
    \item If $\theta\in(1/2,1)$, there exist $\gamma_1>0$ and $\gamma_2>0$ such that
      \begin{equation*}
        |f(\xv_c)-f(\xv_*)| \le \gamma_1 c^{-\frac{1}{2\theta-1}} \qquad\text{ and }\qquad
        \|\xv_c-\xv_*\| \le \gamma_2 c^{-\frac{1-\theta}{2\theta-1}}.
      \end{equation*}
  \end{itemize}
\end{Thm}
\begin{proof}
  Without loss of generality, we assume $\xv_0\in\mathds{B}(\xv_*,\sigma)$. Define
  \begin{equation*}
    \zeta_c:=\sum_{\ell=c}^{\infty} \|\xv_{\ell}-\xv_{\ell+1}\| \ge \|\xv_c-\xv_*\|.
  \end{equation*}
  From Lemma \ref{Lem:5-2} and Theorem \ref{Thm:5-1}, we have
  \begin{eqnarray*}
    \zeta_c &\le& \frac{\max\{1,\overline{\alpha}\}}{\eta C(1-\theta)}|f(\xv_c)-f(\xv_*)|^{1-\theta} \\
      &=& \frac{\max\{1,\overline{\alpha}\}}{\eta C^{1/\theta}(1-\theta)}\left(C|f(\xv_c)-f(\xv_*)|^{\theta}\right)^{(1-\theta)/\theta} \\
      &\le& \frac{\max\{1,\overline{\alpha}\}}{\eta C^{1/\theta}(1-\theta)}\left(\|\gv^1(\xv_c)\|\right)^{(1-\theta)/\theta}.
  \end{eqnarray*}
  It can be obtained from Lemma \ref{Lem:4-2} and \eqref{gc-lowbd} that
  \begin{equation*}
    \|\xv_c-\xv_{c+1}\| = \rho_c\|\gv_c\| \ge \kappa_1\min\{\underline{\alpha},1\}\|\gv^1(\xv_c)\|.
  \end{equation*}
  Hence, by denoting $\kappa_2:=\frac{\max\{1,\overline{\alpha}\}}{\eta C^{1/\theta}(1-\theta)}(\kappa_1\min\{\underline{\alpha},1\})^{-(1-\theta)/\theta}$, we get
  \begin{equation}\label{zeta-1}
    \zeta_c \le \kappa_2\|\xv_c-\xv_{c+1}\|^{(1-\theta)/\theta} = \kappa_2(\zeta_c-\zeta_{c+1})^{(1-\theta)/\theta}.
  \end{equation}
  If $\theta=1/2$, the inequality \eqref{zeta-1} means
  \begin{equation*}
    \zeta_{c+1} \le \frac{\kappa_2-1}{\kappa_2}\zeta_c.
  \end{equation*}
  Hence, the first estimation holds with $\gamma=\zeta_0$ and $\mu=(\kappa_2-1)/\kappa_2$.

  Next, we consider the case $\theta\in(1/2,1)$. Let $\varphi(t):=t^{-\theta/(1-\theta)}$ be a decreasing function for $t>0$. It yields from \eqref{zeta-1} that
  \begin{eqnarray*}
    \kappa_2^{-\theta/(1-\theta)} &\le& \varphi(\zeta_c)(\zeta_c-\zeta_{c+1}) \\
      &=& \int_{\zeta_{c+1}}^{\zeta_c} \varphi(\zeta_c) \mathrm{d}t \\
      &\le& \int_{\zeta_{c+1}}^{\zeta_c} \varphi(t) \mathrm{d}t \\
      &=& \cl{\frac{1-\theta}{1-2\theta}}\left(\zeta_{c+1}^{-(2\theta-1)/(1-\theta)}-\zeta_c^{-(2\theta-1)/(1-\theta)}\right),
  \end{eqnarray*}
  which implies
  \begin{equation*}
    \zeta_{c+1}^{-(2\theta-1)/(1-\theta)}-\zeta_c^{-(2\theta-1)/(1-\theta)}
      \ge \cl{\frac{1-2\theta}{1-\theta}}\kappa_2^{-\theta/(1-\theta)} := \kappa_3>0.
  \end{equation*}
Then, we have
  \begin{equation*}
    \zeta_c^{-(2\theta-1)/(1-\theta)} \ge \kappa_3+\zeta_{c-1}^{-(2\theta-1)/(1-\theta)} \ge\dots
      \ge \kappa_3c+\zeta_0^{-(2\theta-1)/(1-\theta)},
  \end{equation*}
  which means
  \begin{equation*}
    \zeta_c \le \left(\kappa_3c+\zeta_0^{-(2\theta-1)/(1-\theta)}\right)^{-(1-\theta)/(2\theta-1)}
      \le (\kappa_3c)^{-(1-\theta)/(2\theta-1)}.
  \end{equation*}
The last estimation holds by taking $\gamma_2=\kappa_3^{-(1-\theta)/(2\theta-1)}$.

   Let $\xi_c:=|f(\xv_c)-f(\xv_*)|$. From \eqref{fun-lowbd}, \eqref{gc-lowbd}, and Theorem \ref{Thm:5-1}, we have
  \begin{eqnarray*}
    \xi_c-\xi_{c+1} &=& f(\xv_{c+1})-f(\xv_c) \\
      &\ge& \frac{\eta\kappa_1}{\overline{\alpha}}\|\gv_c\|^2 \\
      &\ge& \frac{\eta\kappa_1(\min\{\underline{\alpha},1\})^2}{\overline{\alpha}}\|\gv^1(\xv_c)\|^2 \\
      &\ge& \frac{\eta\kappa_1C^2(\min\{\underline{\alpha},1\})^2}{\overline{\alpha}}\xi_c^{2\theta}
      := \kappa_4 \xi_c^{2\theta}.
  \end{eqnarray*}
Denote $\chi(t):=t^{-2\theta}$ as a decreasing function for $t>0$. It is deduced from the above inequality that
  \begin{equation*}
    \kappa_4 \le \chi(\xi_c)(\xi_c-\xi_{c+1})
      = \int_{\xi_{c+1}}^{\xi_c} \chi(\xi_c) \mathrm{d}t
      \le \int_{\xi_{c+1}}^{\xi_c} \chi(t) \mathrm{d}t
      = \frac{1}{2\theta-1}\left(\xi_{c+1}^{-(2\theta-1)}-\xi_c^{-(2\theta-1)}\right).
  \end{equation*}
  Let $\kappa_5:=(2\theta-1)\kappa_4>0$. We have
  \begin{equation*}
    \xi_c^{-(2\theta-1)} \ge \kappa_5+\xi_{c-1}^{-(2\theta-1)} \ge \dots \ge\kappa_5c+\xi_0^{-(2\theta-1)},
  \end{equation*}
  which means
  \begin{equation*}
    \xi_c \le \left(\kappa_5c+\xi_0^{-(2\theta-1)}\right)^{-1/(2\theta-1)}
      \le (\kappa_5c)^{-1/(2\theta-1)}.
  \end{equation*}
The second estimation is valid if we take $\gamma_1=\kappa_5^{-1/(2\theta-1)}$.
\end{proof}


\section{Numerical experiments}

To evaluate the performance of the proposed algorithm, we implement GPA in MATLAB and use GPA to compute Lagrangians of small and large scale hypergraphs. In our experiments, parameters are set as follows:
$$\eta=0.01, \beta=0.5, \alpha_0=1, \underline{\alpha}=0.001, \text{and} \  \overline{\alpha}=1,000.$$
The algorithm terminates if
$$\|\pi_{\Delta}(\xv_c + \nabla f(\xv_c))-\xv_c\|\le 10^{-8}, |f(\xv_c)-f(\xv_{c-4})|\le 10^{-8},$$
or the number of iteration exceeds one thousand. For each hypergraph, ten random initial points from the legal weighting set $\Delta$ are sampled uniformly. We run the GPA algorithm individually from these starting points, and then choose the best one as our solution. We demonstrate the detailed results in the remainder of this section.

\subsection{A toy example}
First, we consider the toy $3$-graph $G_{toy}$ illustrated in Figure \ref{NE:toy}. The hypergraph $G_{toy}$ has 12 vertices and 13 edges. We compare our method with  function ``fmincon'' in MATLAB optimization tool, which could run the interior point algorithm (IP), sequence quadratic programming (SQP), and the active set method (AS). By solving the optimization, we find that the Lagrangian of $G_{toy}$ is $\lambda(G_{toy})=\frac{1}{16}$ and the associated optimal solution is
\begin{eqnarray*}
  \xv^*_1 &=& (\tfrac{1}{4},\tfrac{1}{4},\tfrac{1}{4},\tfrac{1}{4},0,0,0,0,0,0,0,0)^T, \\
  \xv^*_2 &=& (0,0,0,0,\tfrac{1}{4},\tfrac{1}{4},\tfrac{1}{4},\tfrac{1}{4},0,0,0,0)^T, \\
  \xv^*_3 &=& (0,0,0,0,0,0,0,0,\tfrac{1}{4},\tfrac{1}{4},\tfrac{1}{4},\tfrac{1}{4})^T.
\end{eqnarray*}
It is interesting to see that the Lagrangian always finds the maximal cliques contained in $G_{toy}$. A clique means a complete sub-hypergraph. We will study complete $r$-graph in the next subsection.

\begin{table}[ht]
\caption{CPU time (second).}\label{NE:Toy-Time}
  \centering
   \begin{tabular}{ccccc}
     \toprule
     Methods  & GPA  & IP   & SQP  & AS \\
    \midrule
     time (s) & 0.11 & 3.74 & 0.82 & 1.66\\
   \bottomrule
  \end{tabular}
\end{table}

We report the CPU time of four methods: GPA, IP, SQP, and AS for solving the Lagrangian of $G_{toy}$ in Table \ref{NE:Toy-Time}. It can be seen that GPA is at least seven times faster than fmincon. Since SQP is much faster than IP and AS, we employ fmincon only with SQP in the following experiments.

\subsection{Complete hypergraphs}

It is well-known that the Lagrangian of a complete $r$-graph $K_n^r$ with order $n$ has a closed-form solution
\begin{equation*}
  \lambda(K_n^r)={n \choose r}\frac{1}{n^r}
\end{equation*}
and the associated solution is $\xv^*=(\frac{1}{n},\dots,\frac{1}{n})^T$. In this experiment, we examine complete 3-graphs with order $n$ varying from 10 to 1000 and associated sizes ranging from 120 to 166,167,000.  We compare the results of GPA and SQP  via accuracy of Lagrangian value $|\lambda^*-\hat{\lambda}|$ and the error of optimal solution $\|\xv^*-\hat{\xv}\|_{\infty}$, where $\lambda^*$ and $\xv^*$ are the exact Lagrangian value and the associated optimal solution, and $\hat{\lambda}$ and $\hat{\xv}$ are computed Lagrangian value and computed optimal solution vector.

\begin{table}[H]
\caption{Performance of complete 3-graphs.}\label{NE:completeGraph}
\begin{threeparttable}
\setlength{\tabcolsep}{1mm}{
  \centering
  \begin{tabular}{rrrccrccr}
    \toprule
    \multirow{2}{*}{$n$ } &  \multirow{2}{*}{$m$}  &  \multirow{2}{*}{$\lambda^*$}  &
    \multicolumn{3}{c}{GPA} & \multicolumn{3}{c}{SQP} \\
    \cmidrule(lr){4-6} \cmidrule(lr){7-9}
      &     &   & $|\lambda^*-\hat{\lambda}|$ & $\|\xv^*-\hat{\xv}\|_{\infty}$ & time (s) & $|\lambda^*-\hat{\lambda}|$ & $\|\xv^*-\hat{\xv}\|_{\infty}$ & time (s) \\
   \midrule
       10 &         120 & 0.1200 & $4.0^{-16}$ & $5.1^{-11}$ &     0.14 & $6.1^{-16}$ & $2.4^{-8}$ & 0.33 \\
       18 &         816 & 0.1399 & $1.0^{-15}$ & $8.2^{-11}$ &     0.02 & $3.1^{-16}$ & $1.2^{-8}$ & 0.11 \\
       32 &       4,960 & 0.1514 & $3.7^{-15}$ & $2.2^{-10}$ &     0.04 & $0.0$       & $9.5^{-9}$ & 0.17 \\
       56 &      27,720 & 0.1578 & $8.8^{-15}$ & $8.9^{-11}$ &     0.32 & $5.6^{-17}$ & $1.2^{-8}$ & 0.65 \\
      100 &     161,700 & 0.1617 & $4.2^{-14}$ & $8.4^{-14}$ &     1.39 & $2.5^{-16}$ & $9.1^{-9}$ & 3.07 \\
      178 &     924,176 & 0.1639 & $8.9^{-14}$ & $6.3^{-14}$ &     9.43 & $5.6^{-17}$ & $1.0^{-8}$ & 25.96 \\
      316 &   5,209,260 & 0.1651 & $1.4^{-13}$ & $2.8^{-10}$ &    77.83 & $9.4^{-16}$ & $8.6^{-9}$ & 187.91 \\
      562 &  29,426,320 & 0.1658 & $9.0^{-13}$ & $2.3^{-11}$ &   465.76 & $9.0^{-13}$ & $1.2^{-8}$ & 1,030.73 \\
    1,000 & 166,167,000 & 0.1662 & $1.5^{-12}$ & $2.6^{-8}$  & 5,384.49 & $1.5^{-12}$ & $1.3^{-9}$ & 6,983.94 \\
   \bottomrule
  \end{tabular}
  }
  \end{threeparttable}
\end{table}


Numerical results are illustrated in Table \ref{NE:completeGraph}. First, it is straightforward to see that the value of Lagrangian increases monotonously as the number of vertices of a complete 3-graph enlarges. Second, since GPA is a feasible optimization method, GPA obtains better solution errors than SQP. Finally, GPA  runs faster than SQP when comparing CPU times.

\subsection{Sparse hypergraphs}

\begin{figure}[ht]
  \centering
  \begin{tabular}{ccccc}
    \includegraphics[width=0.2\textwidth]{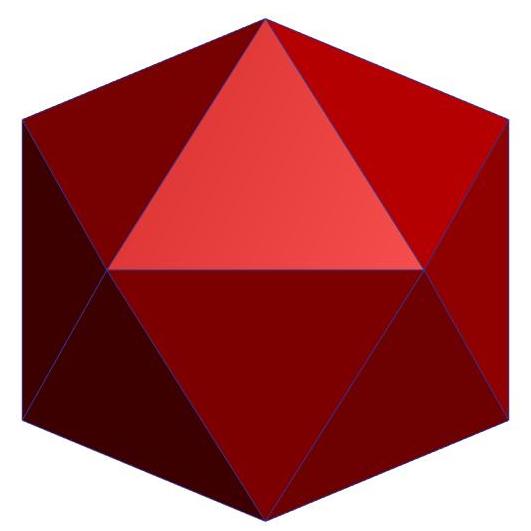} &
    \includegraphics[width=0.2\textwidth]{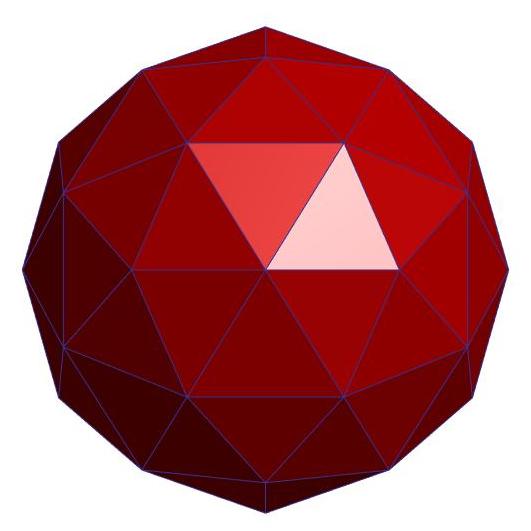} &
    \includegraphics[width=0.2\textwidth]{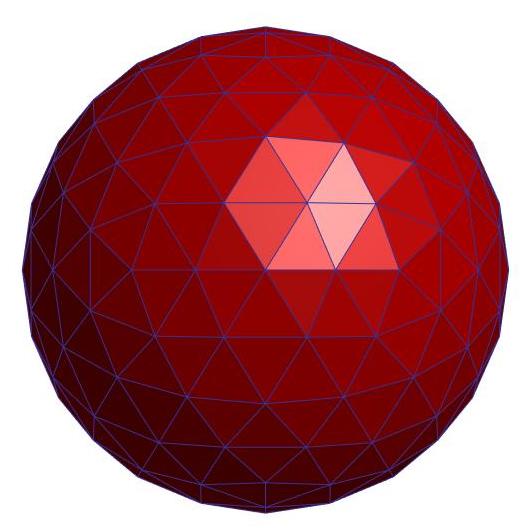} &
    \includegraphics[width=0.2\textwidth]{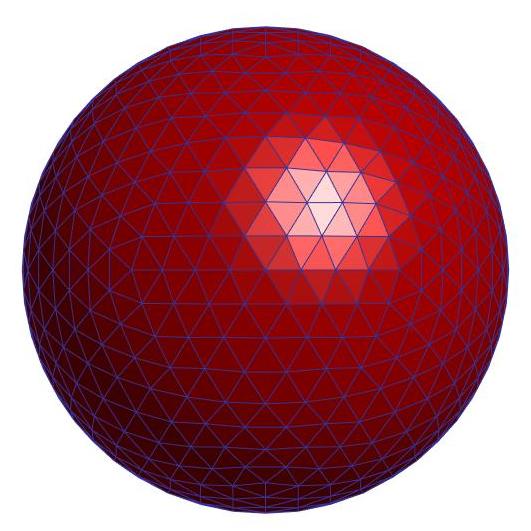} & \raisebox{40pt}{$\dots$} \\
    (a) $\ell=0$ (icosahedron) & (b) $\ell=1$ & (c) $\ell=2$ & (d) $\ell=3$ &
  \end{tabular}
  \caption{Sparse 3-graphs.}\label{NE:triangles}
\end{figure}

In this subsection, we focus on 3-graphs defined on a sphere. As illustrated in Figure \ref{NE:triangles}(a), the icosahedron ($\ell=0$) has 12 vertices and 20 faces. Obviously, each face is a triangle. To approximate the sphere, we subdivide each triangle of the icosahedron into four triangles and obtain the polyhedron in Figure \ref{NE:triangles}(b) with $\ell=1$. By recursively subdividing the triangles, we produce polyhedrons in Figure \ref{NE:triangles}(c) and (d) with $\ell=2$ and $3$, respectively. Then, for each $\ell$, a sparse 3-graph is formed by using the vertex set and the face set of the $\ell$th polyhedron.

\begin{table}[ht]
\caption{Numerical results of sparse 3-graphs.}\label{NE:sparse}
  \centering
  \begin{tabular}{crrcrcr}
     \toprule
    \multirow{2}{*}{$\ell$} &  \multirow{2}{*}{$n$}  &  \multirow{2}{*}{$m$}
                & \multicolumn{2}{c}{GPA}   & \multicolumn{2}{c}{SQP} \\
    \cmidrule(lr){4-5} \cmidrule(lr){6-7}
       & & &        $\hat{\lambda}$ & time (s) & $\hat{\lambda}$ & time (s) \\
    \midrule
        0 &      12 &        20 & 0.037037 &  0.19 & 0.037037 &      0.64 \\
        1 &      42 &        80 & 0.037037 &  0.02 & 0.037037 &      0.14 \\
        2 &     162 &       320 & 0.037037 &  0.03 & 0.037037 &      0.93 \\
        3 &     642 &     1,280 & 0.037037 &  0.09 & 0.037037 &     39.41 \\
        4 &   2,562 &     5,120 & 0.037037 &  0.17 & 0.037037 &    813.32 \\
        5 &  10,242 &    20,480 & 0.037037 &  1.02 & 0.037037 & 64,759.19\\
        6 &  40,962 &    81,920 & 0.037037 &  3.13 & \multicolumn{2}{c}{-- --} \\
        7 & 163,842 &   327,680 & 0.037037 & 16.23 & \multicolumn{2}{c}{-- --} \\
        8 & 655,362 & 1,310,720 & 0.037037 & 62.99 & \multicolumn{2}{c}{-- --} \\
   \bottomrule
  \end{tabular}
  	\noindent{\\\footnotesize{\hspace{-4cm}
  	Here, ``-- --'' means that CPU times exceeds 24 hours.} }
\end{table}

The computation results of Lagrangians  as well as  CPU time are reported in Table \ref{NE:sparse}. No matter how many vertices and edges are involved in the  3-graphs, it seems that all  Lagrangian values are about $\hat{\lambda}=1/27 \approx0.037037$. However, we did not find a solid proof which covers the Lagrangians of sparse hypergraphs illustrated in Figure \ref{NE:triangles}. We note that these sparse hypergraphs are not regular, i.e., degrees of vertices could be five and six, when $\ell\ge1$.
The CPU time shows that our GPA method is thousands times faster than SQP when we solve the sparse 3-graphs with 2,562 vertices. Furthermore, the GPA method is capable of computing the Lagrangian of hypergraphs with millions of edges, which means GPA is powerful for calculating Lagrangians of large scale hypergraphs.

\section{Conclusions}

The Lagrangian of a hypergraph reveals special nature of the adjacency tensor of a uniform hypergraph. In this paper, we designed a gradient projection algorithm for computing the Lagrangian of a uniform hypergraph numerically. Global and local convergence of the proposed algorithm was analyzed. Preliminary numerical experiments illustrated that the proposed algorithm is efficient for large scale hypergraphs.

\begin{acknowledgements}
This work was supported by the National Natural Science Foundation of China (grant No. 11901118 and 62073087), and Suqian Sci$\&$Tech Program (Grant No. Z2020135 and K202112).
\end{acknowledgements}


\noindent\small{\textbf{Data Availability} All data generated or analyzed during this study are included in this manuscript.}





\begin{thebibliography}{plain}


 \bibitem{AtB-09}
 Attouch, H., Bolte, J.:
 On the convergence of the proximal algorithm for nonsmooth functions involving analytic features.
 Math. Program. {\bf 116}(1), 5--16 (2009).
 \url{https://doi.org/10.1007/s10107-007-0133-5}

  \bibitem{ABRS10}
   Attouch, H., Bolte, J., Redont, P., Soubeyran, A.:
   Proximal alternating minimization and projection methods for nonconvex problems: An approach based on the Kurdyka--{\L}ojasiewicz inequality.
   Math. Oper. Res. {\bf 35}(2), 438--457 (2010).
   \url{https://doi.org/10.1287/moor.1100.0449}

  \bibitem{BB88}
  Barzilai, J., Borwein, J.~M.:
  Two-point step size gradientmethods.
  IMA J. Numer. Anal. {\bf 8}, 141--148 (1988).
  \url{https://doi.org/10.1093/imanum/8.1.141}

  \bibitem{BGL16}
  Benson, A.R.,  Gleich, D.F., Leskovec, J.:
  Higher-order organization of complex networks.
  Science {\bf 353}, 163--166 (2016).
  \url{https://doi.org/10.1126/science.aad9029}


  \bibitem{Ber-book}
   Bertsekas, D.P.:
   Nonlinear Programming, 3rd ed.,
   Athena Scientific, Belmont (2016).

  \bibitem{BDL-07}
  Bolte, J.,  Daniilidis, A., Lewis A.:
  The {\L}ojasiewicz inequality for nonsmooth subanalytic functions with applications to subgradient dynamical systems.
  SIAM J. Optim. {\bf 17}(4), 1205--1223 (2007).
  \url{https://doi.org/10.1137/050644641}


  \bibitem{BP13}
  Bul\`{o}, S.R., Pelillo, M.:
  A game-theoretic approach to hypergraph clustering.
  IEEE Trans. Pattern Anal. Mach. Intell. {\bf 35}, 1312--1327 (2013).
  \url{https://doi.org/10.1109/TPAMI.2012.226}



  \bibitem{CCQ16}
  Chang, J.,  Chen, Y., Qi, L.:
  Computing eigenvalues of large scale sparse tensors arising from a hypergraph.
   SIAM J. Sci. Comput. {\bf 38}, A3618--A3643 (2016).
   \url{https://doi.org/10.1109/10.1137/16M1060224}

  \bibitem{CCQY20}
  Chang, J.,  Chen, Y., Qi, L., Yan, H.:
  Hypergraph clustering using a new Laplacian tensor with applications in image processing.
   SIAM J. Imaging Sci. {\bf 13}(3), 1157--1178 (2020).
  \url{https://doi.org/10.1137/19M1291601}

  \bibitem{CDQY18}
  Chang, J.,  Ding, W., Qi, L., Yan, H.:
  Computing the $p$-spectral radii of uniform hypergraphs with applications.
   J. Sci. Comput. {\bf 75}, 1--25 (2018).
  \url{https://doi.org/10.1007/s10915-017-0520-x}

  \bibitem{CQZ17}
   Chen, Y., Qi, L., Zhang, X.:
   The Fiedler vector of a Laplacian tensor for hypergraph partitioning.
   SIAM J. Sci. Comput. {\bf 39}(6), A2508--A2537 (2017).
   \url{https://doi.org/10.1137/16M1094828}

  \bibitem{CY'11}
   Chen, Y., Ye, X.:
   Projection onto a simplex.
   arXiv:1101.6081v2, 1--7 (2011).
   \url{https://doi.org/10.48550/arXiv.1101.6081}

  \bibitem{CoD-12}
   Cooper, J., Dutle, A.:
   Spectra of uniform hypergraphs.
   Linear Algebra Appl. {\bf 436}, 3268--3292 (2012).
   \url{https://doi.org/10.1016/j.laa.2011.11.018}

  \bibitem{FF-89}
   Frankl, P., F\"{u}redi, Z.:
   Extremal problems whose solutions are the blowups of the small witt-designs.
   J. Comb. Theory Ser. A {\bf 52}(1), 129--147 (1989).
   \url{https://doi.org/10.1016/0097-3165(89)90067-8}

  \bibitem{GLM-20}
  Gruslys, V., Letzter,  S., Morrison, N.:
  Hypergraph Lagrangians I: The Frankl-F\"{u}redi conjecture is false.
  Adv. Math. {\bf 365}, 107063 (2020).
   \url{https://doi.org/10.1016/j.aim.2020.107063}

  \bibitem{GLM-21}
  Gruslys, V., Letzter, S., Morrison, N.:
  Lagrangians of hypergraphs II: When colex is best.
  Isr. J. Math.  {\bf 242}, 637--662 (2021).
   \url{https://doi.org/10.1007/s11856-021-2132-2}

  \bibitem{HZ'06}
  Hager, W.W.,  Zhang, H.:
  A new active set algorithm for box constrained optimization.
  SIAM J. Optim. {\bf 17}(2),  526--557 (2006).
  \url{https://doi.org/10.1137/050635225}

  \bibitem{HPY22}
  Hou, J., Pelillo, M., Yuan, H.:
  Hypergraph matching via game-theoretic hypergraph clustering.
  Pattern Recognit. {\bf 125}, 108526 (2022).
  \url{https://doi.org/10.1016/j.patcog.2022.108526}

  \bibitem{HL'18}
  Hu, S., Li, G.:
  Convergence rate analysis for the higher order power method in best rank one approximations of tensors.
  Numer. Math. {\bf 140}, 993--1031 (2018).
  \url{https://doi.org/10.1007/s00211-018-0981-3}

  \bibitem{HPW20}
  Hu, S., Peng, Y., Wu, B.:
  Lagrangian densities of linear forests and Tur\'{a}n numbers of their extensions.
  J. Combin. Des. {\bf 28}, 207--223 (2020).
  \url{https://doi.org/10.1002/jcd.21687}

  \bibitem{HDL21}
   Huang, Y., Dai, Y., Liu, X.:
   Equipping the Barzilai--Borwein method with the two dimensional quadratic termination property.
   SIAM J. Optim. {\bf 31}, 3068--3096 (2021).
   \url{https://doi.org/10.1137/21M1390785}


  \bibitem{K21}
   Koutsoukos, D., Nakandala, S.,  Karanasos, K., Saur, K., Alonso, G., Interlandi, M.:
   Tensors: an abstraction for general data processing.
   Proc. VLDB Endow. {\bf 14}(10), 1797--1804  (2021).
   \url{https://doi.org/10.14778/3467861.3467869}

  \bibitem{Loj-63}
   {\L}ojasiewicz, S.:
   Une propri\'{e}t\'{e} topologique des sous-ensembles analytiques   r\'{e}els.
   Les \'{E}quations aux D\'{e}riv\'{e}es Partielles, 87--89 (1963).


  \bibitem{LZ15}
   Lu, X., Zhang, X.:
   A note on Lagrangians of 4-uniform hypergraphs.
   Ars Combin. {\bf 121}, 329--340 (2015).


  \bibitem{LPL22}
  Luo, X., Peng, J., Liang, J.:
  Directed hypergraph attention network for traffic forecasting.
  IET Intell. Transp. Syst. {\bf 16}, 85--98 (2022).
  \url{https://doi.org/10.1049/itr2.12130}


  \bibitem{MS-65}
  Motzkin, T.S., Straus,  E.G.:
  Maxima for graphs and a new proof of a theorem of Tur\'{a}n.
  Can. J. Math. {\bf 17}, 533--540 (1965).
  \url{https://doi.org/10.4153/cjm-1965-053-6}

  \bibitem{QL-book}
   Qi, L., Luo, Z.:
   Tensor Analysis: Spectral Theory and Special Tensors,
   SIAM, Philadelpia (2017).
   \url{https://doi.org/10.1137/1.9781611974751}


   \bibitem{STZP14}
   Sun, Y,  Tang, Q.,  Zhao, C.,  Peng, Y.:
   On the largest graph-Lagrangian of 3-graphs with fixed number of edges.
   J. Optim. Theory Appl. {\bf 163}, 57--79 (2014).
   \url{https://doi.org/10.1007/s10957-013-0519-x}

  \bibitem{Ta-02}
   Talbot, J.M.:
   Lagrangians of hypergraphs.
   Comb., Probab. Comput. {\bf 11}(2), 199--216 (2022).
   \url{https://doi.org/10.1017/s0963548301005053}

  \bibitem{WC'13}
   Wang, W., Carreira-Perpi\~{n}\'{a}n, M.\'{A}.:
   Projection onto the probability simplex: An efficient algorithm with a simple proof, and an application.
   arXiv:1309.1541v1, 1--5 (2013).
   \url{https://doi.org/10.48550/arXiv.1309.1541}

  \bibitem{ZLZ22}
  Zhang, D., Luo, Y., Yu, Y., Zhao, Q., Zhou, G.:
  Semi-supervised multi-view clustering with dual hypergraph regularized partially shared non-negative matrix factorization.
  Sci. China Technol. Sci. {\bf 65}, 1349--1365 (2022).
   \url{https://doi.org/10.1007/s11431-021-1957-3}

%
%
%
%
%


\end{thebibliography}
\end{document}